\documentclass[10pt,a4paper,reqno]{amsart}
\usepackage[english]{babel}
\usepackage[T1]{fontenc}
\usepackage{amssymb}
\usepackage{amsfonts}
\usepackage{amsmath}
\usepackage{mathptm}

\begin{document}
\thispagestyle{empty}
\pagestyle{plain}
\null{}
\title{\large{On an inequality related to the radial growth of subharmonic functions}}
\author{Juhani Riihentaus}
\date{November 16, 2008}
\maketitle
\begin{center}
{Department of Physics and Mathematics, University of Joensuu\\
P.O. Box 111, FI-80101 Joensuu, Finland \\
juhani.riihentaus@joensuu.fi}
\end{center}

\vspace*{4ex}

\centerline{{{\textbf{ABSTRACT}}}}

\vspace{2ex}

\noindent It is a classical result that every
subharmonic function, defined  and ${\mathcal{ L}}^p$-integrable for some $p$,
$0<p<+\infty$,
on the unit disk $\mathbb{D}$ of the complex plane ${\mathbb{C}}$ is for almost all $\theta$ of the form $o((1-\vert z\vert
)^{-1/p})$, uniformly as
$z\to e^{i\theta}$ in any Stolz domain.
Recently Pavlovi\'c gave a related integral inequality for absolute values of harmonic functions, also defined on the unit disk in the complex plane. 
We generalize Pavlovi\'c's result to so called  quasi-nearly subharmonic functions defined on rather general domains in  ${\mathbb {R}}^n$, $n\geq 2$.

\vspace{3ex}

\noindent



\vspace{8pt}

\noindent{\textbf{Math. Subj. Class.:}} \emph{31B25, 31B05}

\setcounter{equation}{0}

\vspace{3ex}

\thispagestyle{empty}

\noindent{\large{\textbf{1 Introduction}}}

\vspace{3ex}

\noindent{\textbf{1.1 Previous results.}} The following theorem is a special case of  the original  result of Gehring [{4}, 
Theorem~1, p.~77], and 
of Hallenbeck [{5}, Theorems~1 and 2, pp.~117-118],  and of the later and more general results of Stoll 
[{23}, Theorems~1 and 2, pp.~301-302, 307]:

\vspace{1ex}

\noindent\textbf{Theorem A} {\emph{If $u$ is a function harmonic in ${\mathbb{D}}$ such that 
\begin{equation} I(u):=\int_{{\mathbb{D}}}\mid u(z)\mid^p (1-\mid z\mid )^{\beta} \, dm(z)<+\infty ,\end{equation}
where $p>0$, $\beta >-1$, then 
\begin{equation} \lim_{r\rightarrow 1-}\mid u(re^{i\theta })\mid ^p(1-r)^{\beta +1}=0\end{equation}
for almost all $\theta \in [0,2\pi )$.}}

\vspace{1ex}

Observe that  Gehring, Hallenbeck and Stoll in fact considered subharmonic functions and  that the limit in (2) was  uniform in Stolz approach regions 
(in Stoll's result in even more general regions). For a more general result, see [{19}, Theorem, p.~31],  
[{15}, Theorem, p.~233],  [{10}, Theorem~2, p.~73] and [{18}, Theorem~3.4.1, pp.~198-199].  

With the aid of [{{12}}, Theorem~A and Theorem~1, pp.~433-434], Pavlovi\'c showed that the convergence in (2) in Theorem~A is dominated. 
At the same time he pointed out that  whole Theorem~A follows from his result:

\vspace{1ex}

\noindent\textbf{Theorem B} ([12, Theorem~1, pp.~433-434]) {\emph{If $u$ is a function harmonic  in ${\mathbb{D}}$ satisfying}} (1){\emph{, where $p>0$, $\beta >-1$, then 
\begin{equation*} J(u):=\int\limits_{0}^{2\pi}\sup_{0<r<1}\mid u(re^{i\theta })\mid ^p (1-r)^{\beta +1} \, d\theta <+\infty .\end{equation*}
Moreover, there is a constant $C=C_{p,\beta }$ such that $J(u)\leq C\,I(u)$.}}

\vspace{1ex}

The purpose of this note is to point out that, with the aid of  [19, Theorem, p.~31], one can extend Theorem~B considerably: 
Instead of absolute values of harmonic functions on the unit disk ${\mathbb{D}}$ in  the complex plane ${\mathbb{C}}$ we will consider nonnegative 
quasi-nearly subharmonic functions defined on rather general domains of ${\mathbb{R}}^n$, $n\geq 2$. See Theorems~1 and 2 below. 

\vspace{1ex}

First  the necessary  notation and definitions. 

\vspace{2ex}

\noindent{\textbf{1.2 Notation.}}  Our notation is fairly standard, see e.g. [{19}, {21}, {6}]. However, for convenience of the reader 
we recall 
the following.
The common convention $0\cdot \infty =0$ is used. The complex space ${\mathbb{C}}^n$ is identified with the real space  ${\mathbb{R}}^{2n}$, $n\geq 1$.  
In the sequel $D$ is an arbitrary domain in ${\mathbb R}^n$, $n\geq 2$, $D\ne {\mathbb{R}}^n$, whereas $\Omega $ is a bounded domain in ${\mathbb R}^n$ whose boundary 
$\partial \Omega $ is Ahlfors-regular with dimension $d$, $0\leq d\leq n$ (for the definition of this see {\textbf{1.6}} below).
The distance from $x\in D$ to $\partial D$
is denoted by $\delta (x)$. If $\rho >0$ write $D_\rho =\{\,x\in D :\, \delta (x)<\rho \,\}.$
 $B^n(x,r)$ is the Euclidean ball in ${\mathbb{R}}^n$, with center $x$ and radius $r$, and 
$B(x)=B^n(x,\frac{1}{3}\delta (x))$. We write $B^n=B(0,1)$ and $S^{n-1}=\partial B^n$.  $m$ is the Lebesgue measure in ${\mathbb{R}}^n$, 
and $\nu _n=m(B^n)$.
${\mathcal L}^1_{\mathrm{loc}}(D)$ is the space of
locally (Lebesgue) integrable functions on $D$. The \mbox{$d$-dimensional} Hausdorff (outer) measure in ${\mathbb R}^n$ is denoted by $H^d$,
$0\leq d\leq n$. 
Our constants $C$ and $K$ are always positive, mostly $\geq 1$ and they may vary from line to line. (One exception: In the proof of Theorem~2 we 
 write $K$ for $\partial \Omega $, just in order to follow our previous notation in [{{19}}].) On the other hand, $C_0$ and $r_0$ are fixed 
constants which  are involved with the used (and thus   fixed)  admissible function $\varphi $ (see {\textbf{1.5}} (5) below). Similarly, if $\alpha >0$ is given,
$C_1=C_1(C_0,\alpha )$, $C_2=C_2(C_0,\alpha )$ and  $C_3=C_3(C_0,\alpha )$ are fixed constants, coming directly from 
[{19}, Lemma~2.3, pp.~32-33] or [{15}, Lemma~2.3, p.~234],  and thus defined already there.

\vspace{2ex}

\noindent{\textbf{1.3 Nearly subharmonic functions.}} We recall that an upper semicontinuous function $u:\, D\rightarrow [-\infty ,+\infty )$ is \emph{subharmonic} if 
for all $\overline{B^n(x,r)}\subset D$,
\[u(x)\leq \frac{1}{\nu _n\, r^n}\int\limits_{B^n(x,r)}u(y)\, dm(y).\]
The function $u\equiv -\infty $  is considered  subharmonic. 

We say that a function 
$u:\, D\rightarrow [-\infty ,+\infty )$ is \emph{nearly subharmonic}, if $u$ is Lebesgue measurable, $u^+\in {\mathcal{L}}^1_{\textrm{loc}}(D)$, 
and for all $\overline{B^n(x,r)}\subset D$,  
\begin{equation*}u(x)\leq \frac{1}{\nu _n\, r^n}\int\limits_{B^n(x,r)}u(y)\, dm(y).\end{equation*}
Observe that in the standard definition of nearly subharmonic functions one uses the slightly stronger  assumption that 
$u\in {\mathcal{L}}^1_{\textrm{loc}}(D)$, see e.g. [{6}, p.~14]. However, our above, slightly 
more general definition seems to be  more useful, see   [{21}, Proposition~2.1~(iii) 
and Proposition~2.2~(vi), (vii), pp.~54-55].

\vspace{2ex}

\noindent{\textbf{1.4 Quasi-nearly subharmonic functions.}}  
A Lebesgue measurable function $u:\,D \rightarrow 
[-\infty ,+\infty )$ is \emph{$K$-quasi-nearly subharmonic}, if  $u^+\in{\mathcal{L}}^{1}_{\textrm{loc}}(D)$ and if there is a 
constant $K=K(n,u,D)\geq 1$
such that for all   $\overline{B^n(x,r)}\subset D$,    
\begin{equation} u_M(x)\leq \frac{K}{\nu _n\,r^n}\int\limits_{B^n(x,r)}u_M(y)\, dm(y)\end{equation}
for all $M\geq 0$, where $u_M:=\sup\{u,-M\}+M$. A function $u:\, D\rightarrow [-\infty ,+\infty )$ is \emph{quasi-nearly subharmonic}, if $u$ is 
$K$-quasi-nearly subharmonic for some $K\geq 1$.

\vspace{1ex}

A Lebesgue  measurable function 
$u:\,D \rightarrow [-\infty ,+\infty )$ is \emph{$K$-quasi-nearly subharmonic n.s. (in the narrow sense)}, if $u^+\in{\mathcal{L}}^{1}_{\textrm{loc}}(D)$ and if there is a 
constant $K=K(n,u,D)\geq 1$
such that for all $\overline{B^n(x,r)}\subset D$,   
\begin{equation} u(x)\leq \frac{K}{\nu _n\,r^n}\int\limits_{B^n(x,r)}u(y)\, dm(y).\end{equation}
A function $u:\, D\rightarrow [-\infty ,+\infty )$ is \emph{quasi-nearly subharmonic n.s.}, if $u$ is 
$K$-quasi-nearly subharmonic n.s.  for some $K\geq 1$. 

\vspace{1ex}

Quasi-nearly subharmonic functions (perhaps with a different terminology), or,  essentially, perhaps just functions satisfying a certain  generalized mean value 
inequality, more or less  of the form  (3) or (4) above, have previously 
been considered or used   at least in [{3}, {25}, {8}, {14}, {24}, {5}, {11}, {9}, {23}, {15}, {10}, 
{16}, {17}, {18}, {13}, {19}, {20}, {21}, {7}]. 
We recall here only that this  function class 
includes, among \mbox{others,}  subharmonic functions, and, more generally,  quasisubharmonic and nearly subharmonic functions (for the definitions 
of these,    see above and e.g. [{6}]),
also functions satisfying certain natural  growth conditions, especially  
certain eigenfunctions,  polyharmonic functions, subsolutions of certain general elliptic equations.   Also, the class of Harnack functions is included, 
thus, among others, nonnegative harmonic functions 
as well as nonnegative solutions of some elliptic equations. In particular, the partial differential equations associated with quasiregular mappings 
belong to this family of elliptic equations, see Vuorinen [{26}].  Observe that already Domar  [{2}] has pointed out the relevance of the 
class of (nonnegative) quasi-nearly subharmonic functions.  

\vspace{1ex}

To motivate the reader still further, we recall here  the following, see e.g. [{13}, Proposition~1, Theorem~A, Theorem~B, p.~91] and 
[{21}, Proposition~2.1 and Proposition~2.2, pp.~54-55]:  

\begin{itemize}
\item[(i)] \emph{A $K$-quasi-nearly subharmonic function n.s. is $K$-quasi-nearly subharmonic, but not necessarily conversely.}
\item[(ii)] \emph{A nonnegative Lebesgue measurable function is $K$-quasi-nearly subharmonic if and only if it is $K$-quasi-nearly subharmonic n.s.} 
\item[(iii)] \emph{A Lebesgue measurable function is $1$-quasi-nearly subharmonic if and only if it is $1$-quasi-nearly subharmonic n.s. and 
if and only if it is nearly subharmonic (in the sense  defined above).}
\item[(iv)] \emph{If  $u:\, D\rightarrow [0 ,+\infty )$   is quasi-nearly subharmonic and $p>0$, then  $u^p$ is quasi-nearly subharmonic. Especially, if 
$h: \,D\rightarrow {\mathbb{R}}$ is harmonic and $p>0$, then $\mid h\mid ^p$ is quasi-nearly subharmonic.}
\item[(v)] \emph{If  $u:\, D\rightarrow [-\infty ,+\infty )$   is quasi-nearly subharmonic n.s., then either $u\equiv -\infty $ or $u$ is finite almost 
everywhere in $D$, and} 
$u\in {\mathcal{L}}^1_{\textrm{loc}}(D)$.
\end{itemize}

\vspace{2ex}

\noindent{\textbf{1.5 Admissible functions.}} A function
$\varphi :[0,+\infty )\rightarrow [0,+\infty )$ is {\emph{admissible}}, if
it is strictly increasing, surjective and  there are constants
$C_0=C_0(\varphi )\geq 1$ and
$r_0>0$ such that
\begin{equation} \varphi
(2t)\leq C_0\, \varphi
(t)
\, \, \, \, 
{\textrm{ {and}}}\, \, \, \, \, \,  \varphi ^{-1} (2s)\leq C_0\, \varphi ^{-1}(s)\end{equation}
for all
$s,\, t$, $0\leq s,\, t\leq r_0$.

\vspace{1ex}

Functions $\varphi _1(t)=t^\tau $, $\tau > 0$, or, more generally, nonnegative,  increasing surjective functions $\varphi_2(t)$ which satisfy
 the $\varDelta_2$-condition and for which the functions
$t\mapsto \frac{\varphi_2(t)}{t}$ are increasing,
are examples of admissible functions.
Further examples are
$\varphi_3(t)=c\, t^{\alpha}[\log (\delta +t^{\gamma})]^{\beta},$
where
$c>0$,
$\alpha
>0$,
$\delta
\geq
1$,  and
$\beta
,\gamma
\in
{\mathbb
{R}}$  are such that
$\alpha +\beta \gamma >0$. For more examples, see [{15}, {18}].

\vspace{1ex}

Let $\varphi :[0,+\infty )\rightarrow [0,+\infty )$ be an
admissible function and let $\alpha >0$. One says that $\zeta \in \partial D$ is
$(\varphi ,\alpha )$-{\it {accessible}}, shortly \emph{accessible}, if
\[ \Gamma_{\varphi}(\zeta ,\alpha )\cap B^n(\zeta ,\rho )\ne
\emptyset \] for all
$\rho
 >0$. Here
\[ \Gamma_{\varphi} (\zeta ,\alpha )=\{\, x\in D \, :\, \varphi (\vert
x-\zeta \vert )<\alpha \, \delta (x)\, \},\]
and it is called a  $(\varphi ,\alpha )$-{\textit {approach region}}, shortly an \emph{approach region}, \emph{in}
$D$ {\textit {at}} $\zeta$. Choosing $\varphi (t)=t$ (in the case of the unit disk ${\mathbb{D}}$ of the complex plane ${\mathbb{C}}$)  one gets 
the familiar Stolz approach region. Choosing $\varphi (t)=t^\tau $, $\tau \geq 1$, say, one gets more general approach regions, see [{23}].

\vspace{2ex}

\noindent{\textbf{1.6}} Let $0\leq d\leq n$. A set $E\subset {\mathbb{R}}^n$ is \emph{Ahlfors-regular with dimension $d$} if it is closed and there is 
a constant $C_4>0$ so that 
\[ C_4^{-1}r^d\leq H^d(E\cap B^n(x,r))\leq C_4r^d\]
for all $x\in E$ and $r>0$. The smallest constant $C_4$ is called the \emph{regularity constant} for $E$. Simple examples of Ahlfors-regular sets 
include $d$-planes and $d$-dimensional Lipschitz graphs. Also certain Cantor sets and self-similar sets are Ahlfors-regular. For more details, 
see [{1}, pp.~9-10].

\vspace{4ex}

\noindent{\large\textbf{2 The results}}

\vspace{3ex}

\noindent\textbf{2.1} First  a partial generalization to Pavlovi\'c's result [{12}, Theorem~1, pp.~433-434] or  Theorem~B above. Observe that though the constant 
$C$ below in (6) does depend on $K$, it is, nevertheless, otherwise independent of the \mbox{($K$-)}quasi-nearly subharmonic \mbox{function $u$.}

\vspace{1ex}

\noindent\textbf{Theorem~1} {\emph{Let $\Omega $  be a  domain in ${\mathbb{R}}^n$, $n\geq 2$, $\Omega \ne {\mathbb{R}}^n$,  such that its boundary $\partial \Omega$
is Ahlfors-regular with dimension $d$,   $0\leq d\leq n$. Let $u:\Omega \rightarrow [0,+\infty )$ be a $K$-quasi-nearly subharmonic function. Let 
$\varphi :[0,+\infty )\rightarrow [0,+\infty )$ be an admissible function, with constants $r_0$ and $C_0$. Let $\alpha >0$ be arbitrary. Let 
$\rho _0:=\min \{\, r_0/2^{1+\alpha },r_0/2^{3\alpha C_0},\varphi (r_0)/\alpha \,\}$. Let $\gamma \in {\mathbb{R}}$ be such that 
\begin{equation*} \int_{\Omega}\delta (x)^{\gamma}\,u(x)\, dm(x)<+\infty . \end{equation*}
Then there is a constant $C=C(n,\Omega ,d, \varphi , \alpha ,\gamma ,K)$  such that for all $\rho \leq \rho _0$,
\begin{equation*} \int\limits_{\partial \Omega }
\sup_{x\in \Gamma _{\varphi ,\rho }(\zeta ,\alpha )}\{\, \delta (x)^{n+\gamma }[\varphi ^{-1}(\delta (x))]^{-d}u(x)\,\}
\,dH^d(\zeta )\leq C\int\limits_{\Omega_{\rho'}  }\delta (x)^\gamma\,u(x)\,dm(x),
\end{equation*}
where $\rho '=\frac{4}{3}\rho $ and 
\[\Gamma_{\varphi ,\rho } (\zeta ,\alpha )=\{\, x\in \Gamma_{\varphi } (\zeta ,\alpha ) :\, \delta (x)<\rho \, \}.\]
}}

\vspace{1ex}

\noindent\textbf{Proof.}  Proceeding as in 
[{19}, proof of Theorem (with $\psi =id$), pp.~31-35] (cf. 
[{15},  proof of  Theorem, pp.~235-237]) and choosing $K=\partial \Omega $, one obtains 
\begin{equation*} \int\limits_{\partial \Omega }M_{\rho}^{\partial \Omega } (\zeta )\,dH^d(\zeta )\leq C\int\limits_{\Omega_{\rho'}}
\delta (x)^\gamma\,u(x)\,dm(x)
\end{equation*}
where $\rho '=\frac{4}{3}\rho $ and $M_\rho ^{\partial \Omega }:\partial \Omega \rightarrow [0,+\infty ]$,
\begin{equation*}M^{\partial \Omega}_{\rho}(\zeta)=
\sup_{x\in \Gamma_{\varphi ,\rho }(\zeta ,\alpha )}\frac{\delta (x)^{n+\gamma}u(x)}{[\varphi ^{-1}(\delta
(x))]^{d}+H^d(B^n(x,C_1C_2\, \varphi^{-1}(\delta (x)))\cap \partial \Omega )}.\end{equation*}

Here and below the constants $C_1=C_1(C_0,\alpha )$, $C_2=C_2(C_0,\alpha )$ and $C_3=C_3(C_0,\alpha )$ are, as pointed out above,  directly from 
[{19}, proof of Lemma~2.3, pp.~32-33] or [{15}, proof of Lemma~2.3, pp.~234-235]. By this lemma  one has, for each $\zeta \in \partial \Omega $ and for each 
$x\in \Gamma _{\varphi ,\rho }(\zeta ,\alpha )$, $B^n(x,C_1C_2\varphi ^{-1}(\delta (x)))\subset B^n(\zeta ,C_1C_2C_3\varphi ^{-1}(\delta (x)))$.
Since $\partial \Omega $ is Ahlfors-regular with dimension $d$, we have
\begin{equation*}H^d(B^n(\zeta ,C_1C_2C_3\varphi ^{-1}(\delta (x)))\cap \partial \Omega )\leq C_4[C_1C_2C_3\varphi ^{-1}(\delta (x))]^d\end{equation*}
where also $C_4$ is a fixed constant. Therefore
\begin{equation*}\begin{split}M^{\partial \Omega}_{\rho}(\zeta)&=
\sup_{x\in \Gamma_{\varphi ,\rho }(\zeta ,\alpha )}\frac{\delta (x)^{n+\gamma}u(x)}{[\varphi ^{-1}(\delta
(x))]^{d}+H^d(B^n(x,C_1C_2\, \varphi^{-1}(\delta (x)))\cap \partial \Omega )}\\
&\geq \sup_{x\in \Gamma_{\varphi ,\rho }(\zeta ,\alpha )}\frac{\delta (x)^{n+\gamma}u(x)}{[\varphi ^{-1}(\delta
(x))]^{d}+H^d(B^n(\zeta ,C_1C_2C_3\, \varphi^{-1}(\delta (x)))\cap \partial \Omega )}\\
&\geq \sup_{x\in \Gamma_{\varphi ,\rho }(\zeta ,\alpha )}\frac{\delta (x)^{n+\gamma}u(x)}{[\varphi ^{-1}(\delta
(x))]^{d}+C_4(C_1C_2C_3)^d[\varphi^{-1}(\delta (x))]^d} \\
&\geq \frac{1}{1+(C_1C_2C_3)^dC_4}
\sup_{x\in \Gamma _{\varphi ,\rho }(\zeta ,\alpha )}\{\, \delta (x)^{n+\gamma }[\varphi ^{-1}(\delta (x))]^{-d}u(x)\,\}.
\end{split}\end{equation*}
Hence
\begin{equation*} \int\limits_{\partial \Omega }
\sup_{x\in \Gamma _{\varphi ,\rho }(\zeta ,\alpha )}\{\, \delta (x)^{n+\gamma }[\varphi ^{-1}(\delta (x))]^{-d}u(x)\,\}
\,dH^d(\zeta )\leq C\int\limits_{\Omega_{\rho'}}
\delta (x)^\gamma\,u(x)\,dm(x),
\end{equation*}
concluding the proof. 
\hfill\qed

\vspace{2ex}

\noindent\textbf{2.2} Theorem~1 seems to be useful in many situations. For example, with the aid of it one gets the following improvements to 
Pavlovi\'c's result [{12}, Theorem~1, pp.~433-434] or  Theorem~B above:

\vspace{1ex}
    
\noindent\textbf{Theorem~2} {\emph{Let $\Omega $, $d$, $u$, $\varphi $, $\alpha $,  $\gamma $ and $\rho _0$ be as above in Theorem~1. Suppose 
moreover that $H^d(\partial \Omega )<+\infty $.
Then there is a constant $C=C(n,\Omega ,d, \varphi , \alpha ,\gamma ,K)$  such that 
\begin{equation*} \int\limits_{\partial \Omega }\sup_{x\in \Gamma _\varphi (\zeta ,\alpha )}
\{\delta (x)^{n+\gamma }[\varphi ^{-1}(\delta (x))]^{-d}u(x)\}\,dH^d(\zeta )\leq C\int\limits_{\Omega}\delta (x)^\gamma\,u(x)\,dm(x).
\end{equation*}
}}
\vspace{1ex}

\noindent\textbf{Proof.} By Theorem~1 (we may clearly assume that $\int_{\Omega}\delta (x)^{\gamma}\,u(x)\, dm(x)<+\infty $),
\begin{equation*} \int\limits_{\partial \Omega }
\sup_{x\in \Gamma _{\varphi ,\rho_0 }(\zeta ,\alpha )}\{\, \delta (x)^{n+\gamma }
[\varphi ^{-1}(\delta (x))]^{-d}u(x)\,\}\,dH^d(\zeta )\leq C\int\limits_{\Omega_{\rho_0 '} }\delta (x)^\gamma\,u(x)\,dm(x).
\end{equation*}

Write
\[ \Gamma^c_{\varphi ,\rho _0} (\zeta ,\alpha ):=\{\, x\in \Gamma _{\varphi }(\zeta ,\alpha ):\, \, \delta (x)\geq \rho _0\}.\]
Since  
\begin{equation*}
\begin{split}
\sup_{x\in \Gamma _{\varphi  }(\zeta ,\alpha )}\{\, \delta (x)^{n+\gamma }[\varphi ^{-1}(\delta (x))]^{-d}u(x)\,\}
&
\leq \sup_{x\in \Gamma ^c_{\varphi ,{\rho}_0 }(\zeta ,\alpha )}\{\, \delta (x)^{n+\gamma }\varphi ^{-1}(\delta (x))]^{-d}u(x)\,\}\\
&
+\sup_{x\in \Gamma _{\varphi ,{\rho}_0}(\zeta ,\alpha )}\{\, \delta (x)^{n+\gamma }\varphi ^{-1}(\delta (x))]^{-d}u(x)\,\},
\end{split}
\end{equation*}
we obtain:
\begin{equation*}\begin{split}&\int\limits_{\partial \Omega }\sup_{x\in \Gamma _{\varphi  }(\zeta ,\alpha )}
\{\, \delta (x)^{n+\gamma }[\varphi ^{-1}(\delta (x))]^{-d}u(x)\,\}\, dH^d(\zeta )\\
&\leq \int\limits_{\partial \Omega }\sup_{x\in \Gamma ^c_{\varphi ,{\rho}_0 }(\zeta ,\alpha )}
\{\, \delta (x)^{n+\gamma }[\varphi ^{-1}(\delta (x))]^{-d}u(x)\,\}\, dH^d(\zeta )\\
&+\int\limits_{\partial \Omega }\sup_{x\in \Gamma _{\varphi ,{\rho}_0}(\zeta ,\alpha )}
\{\, \delta (x)^{n+\gamma }[\varphi ^{-1}(\delta (x))]^{-d}u(x)\,\}\,dH^d(\zeta )\\
&\leq \int\limits_{\partial \Omega }\sup_{x\in \Gamma ^c_{\varphi ,{\rho}_0 }(\zeta ,\alpha )}
\{\, \delta (x)^{n+\gamma }[\varphi ^{-1}(\delta (x))]^{-d}u(x)\,\}\, dH^d(\zeta )+
C\int\limits_{\Omega_{\rho_0 '} } \delta (x)^{\gamma }u(x)\,dm(x)\\
&\leq \int\limits_{\partial \Omega }\sup_{x\in \Gamma ^c_{\varphi ,{\rho}_0 }(\zeta ,\alpha )}
\{\, \delta (x)^{n+\gamma }[\varphi ^{-1}(\delta (x))]^{-d}u(x)\,\}\, dH^d(\zeta )+
C\int\limits_{\Omega } \delta (x)^{\gamma }u(x)\,dm(x).
\end{split}\end{equation*}
It remains to show that  
\begin{equation*} \int\limits_{\partial \Omega }\sup_{x\in \Gamma^c _{\varphi ,{\rho}_0}(\zeta ,\alpha )}
\{\delta (x)^{n+\gamma }[\varphi ^{-1}(\delta (x))]^{-d}u(x)\}\,dH^d(\zeta )\leq C\int\limits_{\Omega}\delta (x)^\gamma\,u(x)\,dm(x)
\end{equation*}
for some $C=C(n,\Omega ,d,\varphi ,\alpha ,\gamma ,K)$.
For all $x\in  \Gamma^c _{\varphi ,\rho _0 }(\zeta ,\alpha )$ we have 
\[u(x)\leq \frac{K}{\nu _n(\frac{\delta (x)}{3})^n}\int\limits_{B(x)}u(y)\, dm(y).\]
Using also the facts that  $\frac{2}{3}\delta (x)\leq \delta (y)\leq \frac{4}{3}\delta (x)$
for all $y\in B(x)$, one gets easily:
\begin{equation*}\begin{split} \int\limits_{\partial \Omega }&\sup_{x\in \Gamma^c _{\varphi ,{\rho}_0}(\zeta ,\alpha )}
\{\delta (x)^{n+\gamma }[\varphi ^{-1}(\delta (x))]^{-d}u(x)\}\,dH^d(\zeta )\\
&\leq  \int\limits_{\partial \Omega }\sup_{x\in \Gamma^c _{\varphi ,{\rho}_0}(\zeta ,\alpha )}
\{\delta (x)^{n+\gamma }[\varphi ^{-1}(\delta (x))]^{-d}\frac{K}{\nu _n(\frac{\delta (x)}{3})^n}\int\limits_{B(x)}u(y)\, dm(y)\}dH^d(\zeta )\\
&\leq \frac{3^nK}{\nu _n} \int\limits_{\partial \Omega }\sup_{x\in \Gamma^c _{\varphi ,{\rho}_0}(\zeta ,\alpha )}
\{\delta (x)^{\gamma }[\varphi ^{-1}(\delta (x))]^{-d}\int\limits_{B(x)}u(y)\, dm(y)\}dH^d(\zeta )\\
&\leq \left( \frac{3}{2}\right)^{\mid \gamma \mid }\frac{3^nK}{\nu _n} \int\limits_{\partial \Omega }\sup_{x\in \Gamma^c _{\varphi ,{\rho}_0}(\zeta ,\alpha )}
\{[\varphi ^{-1}(\delta (x))]^{-d}\int\limits_{B(x)}\delta (y)^\gamma u(y)\, dm(y)\}dH^d(\zeta )\\
&\leq \frac{3^{\mid \gamma \mid +n}K}{2^{\mid \gamma \mid }\nu _n}[\varphi ^{-1}(\rho _0)]^{-d}H^d(\partial \Omega )
 \int\limits_{\Omega }\delta (y)^\gamma u(y)\, dm(y).
\end{split}\end{equation*}
Thus
\begin{equation*} \int\limits_{\partial \Omega }
\sup_{x\in \Gamma _{\varphi }(\zeta ,\alpha )}\{\, \delta (x)^{n+\gamma }
[\varphi ^{-1}(\delta (x))]^{-d}u(x)\,\}\,dH^d(\zeta )\leq C\int\limits_{\Omega}\delta (x)^\gamma\,u(x)\,dm(x),
\end{equation*}
concluding the proof.
\null{}\quad \hfill \qed

\vspace{1ex}

\noindent\textbf{Corollary} {\emph{Let $u: B^n\rightarrow [0,+\infty )$ be a subharmonic function 
and let $p>0$, $\alpha >1$ and  $\gamma > -1-\max \{\, (n-1)(1-p),0\,\}$. Then there is a constant $C=C(n,\gamma ,p,\alpha )$ such that
\begin{equation*} \int\limits_{S^{n-1}}\sup_{x\in \Gamma _{id}(\zeta ,\alpha )}\{(1-\mid x\mid )^{\gamma +1}u(x)^p\}\,d\sigma (\zeta )
\leq C\int\limits_{B^n}(1-\mid x\mid )^\gamma u(x)^p\,dm(x).
\end{equation*}
Here $id$ is the identity mapping of ${\mathbb{R}}^n$ and $\sigma $ is the spherical (Lebesgue) measure in $S^{n-1}$.
}}
\vspace{1ex}

\noindent\textbf{Remark} Observe that Suzuki [{{24}}, Theorem~2, pp.~272-273] has shown the following: If $p>0$ and 
$\gamma \leq -1-\max \{\, (n-1)(1-p),0\,\}$,  then
the only nonnegative  subharmonic function  on a bounded domain $D$ of ${\mathbb{R}}^n$ with ${\mathcal{C}}^2$ boundary satisfying 
\begin{equation} \int_{D}\delta (x)^{\gamma}\,u(x)^p\, dm(x)<+\infty  \end{equation}
is the zero function. 
On the other hand, if $p>0$ and $\gamma >-1-\max \{\, (n-1)(1-p),0\,\}$, then  there exist nonnegative non-zero subharmonic functions
on $D=B^n$ satisfying (6). 

\vspace{4ex}

\noindent{\large{\textbf{References}}}

\begin{enumerate}

\item[{[1]}] 
G. David and S. Semmes,  \emph{Analysis of and on Uniformly Rectifiable Sets},  
Math. Surveys and Monographs~38,  Amer.  Math.  Soc., Providence, Rhode Island (1991). \par

\item[{[2]}]
 Y. Domar, \emph{On the existence of a largest subharmonic minorant of a given function}, 
Ark. mat. {\textbf{3}}, nr. 39 (1957), 429--440.\par

\item[{[3]}]
 C.~Fefferman and E.M.~Stein, H$^p$ \emph{spaces of several variables}, 
Acta Math. \textbf{129} (1972),   \mbox{137--192.}\par

\item[{[4]}]
 F.W. Gehring, \emph{On the radial order of subharmonic functions}, 
J. Math. Soc. Japan \textbf{9} (1957), \mbox{77--79.}\par

\item[{[5]}]
D.J. Hallenbeck, \emph{Radial growth of subharmonic functions}, in:
Pitman Research Notes  262 (1992), \mbox{113--121.}\par

\item[{[6]}] 
 M. Hervé, \emph{Analytic and Plurisubharmonic Functions in Finite and Infinite
Dimensional Spaces}, Lecture Notes in Math. 198,  Springer, Berlin $\cdot$ Heidelberg $\cdot$ New York (1971).\par

\item[{[7]}]
 V.~Koji\'c, \emph{Quasi-nearly subharmonic functions and conformal mappings}, 
Filomat. {\textbf{21}}, no. 2 (2007), 243--249.\par
 
\item[{[8]}] \"U.~Kuran, \emph{Subharmonic behavior of $\mid h\mid ^p$, ($p>0$, $h$ harmonic)}, 
J. London Math. Soc. (2) \textbf{8} (1974),  \mbox{529--538.}\par

\item[{[9]}]
  Y.~Mizuta, \emph{Potential Theory in Euclidean Spaces},
 Gaguto International Series, Mathematical Sciences and Applications 6, Gakk$\bar{{\textrm{o}}}$tosho Co., Tokyo (1996).\par

\item[{[10]}]
 Y.~Mizuta, \emph{Boundary limits of functions in weighted Lebesgue or Sobolev classes}, 
Revue Roum. Math. Pures  Appl. \textbf{46} (2001),  \mbox{67--75.}\par

\item[{[11]}]
 M. Pavlovi\'c, \emph{On subharmonic behavior and oscillation of functions
 on balls in ${\mathbb R}^n$}, Publ. de l'Inst. Math., Nouv. sér. \textbf{55(69)} (1994),  \mbox{18--22.}\par

 \item[{[12]}] M. Pavlovi\'c, \emph{An inequality related to the  Gehring-Hallenbeck theorem on radial limits of functions in  harmonic Bergman spaces}, 
Glasgow Math. J. {\textbf{50}}, no.~3 (2008), \mbox{433-435}

\item[{[13]}]
 M.~Pavlovi\'c and J.~Riihentaus, \emph{Classes of quasi-nearly subharmonic functions}, Potential Anal. {\textbf{29}} (2008), 
\mbox{89--104.}\par

\item[{[14]}]
 J. Riihentaus, \emph{On a theorem of Avanissian--Arsove},
Expo.  Math. \textbf{7} (1989),  \mbox{69--72.}\par

\item[{[15]}]
J. Riihentaus, \emph{Subharmonic functions: non-tangential and 
tangential boundary
behavior}, in:   Function Spaces, Differential Operators and Nonlinear Analysis \mbox{(FSDONA'99)}, Proceedings of the Sy\"ote Conference 1999, 
  V.~Mustonen, J.~R\'akosnik (eds.), 
 Math. Inst., Czech Acad. Science,  Praha (2000), \mbox{229--238.} \newline\mbox{(ISBN 80-85823-42-X)}\par

\item[{[16]}]
J.~Riihentaus, \emph{A generalized mean value inequality for 
subharmonic functions}, Expo. Math. \textbf{19} (2001), \mbox{187--190.}\par

\item  [{[17]}]
 J.~Riihentaus, \emph{Subharmonic functions, mean value inequality, boundary behavior, nonintegrability and exceptional sets}, 
in: International Workshop on Potential Theory and Free Boundary Flows, Kiev, Ukraine, August 19-27, 2003, 
Trans. Inst. Math. Nat. Acad.  Sci. Ukr. {\textbf{1}}, no. {\textbf{1}} (2004), 169--191.\par
 
\item[{[18]}]
J.~Riihentaus, \emph{Weighted boundary behavior and nonintegrability of subharmonic functions}, in: International Conference on 
Education and Information Systems: Technologies and Applications (EISTA'04), Orlando, Florida, USA, July 21-25, 2004, Proceedings,
M.~Chang, Y-T.~Hsia, F.~Malpica, M.~Suarez, A.~Tremante, F.~Welsch (eds.), Vol.~II (2004), \mbox{pp. 196--202.} \mbox{(ISBN 980-6560-11-6)}\par

\item[{[19]}]
J.~Riihentaus, \emph{A weighted boundary limit result for subharmonic functions}, Adv. Algebra and Analysis {\textbf{1}} (2006), 
\mbox{27--38.}\par

\item[{[20]}]
J.~Riihentaus, \emph{Separately quasi-nearly subharmonic functions}, in:  Complex Analysis and Potential Theory, Proceedings of the 
Conference Satellite to ICM~2006, Tahir Aliyev Azero$\breve{\textrm{g}}$lu, Promarz M. Tamrazov (eds.), Gebze Institute of Technology, Gebze, Turkey, 
September 8-14,
 2006, World Scientific, Singapore (2007), 156--165.\par

\item[{[21]}]
 J.~Riihentaus, \emph{Subharmonic functions, generalizations and separately subharmonic functions}, 
XIV-th Conference on Analytic Functions, July 22-28, 2007, Che\l m, Poland, in: Sci. Bull. Che\l m, Sect. Math. Comp. Sci. {\textbf{2}} (2007), 49--76. 
\newline(ISBN 978-83-61149-24-8) (arXiv:math/0610259v5 [math.AP] 8 Oct 2008)\par
 
\item[{[22]}]
J.~Riihentaus, \emph{Quasi-nearly subharmonicity and separately quasi-nearly subharmonic functions}, J. Inequal. Appl. {\textbf{2008}} (2008), Article ID 149712, 15~pages \newline(doi: 10.1155/2008/149712). 

\item[{[23]}]
M.~Stoll, \emph{Weighted tangential boundary limits of subharmonic functions on
domains in ${\mathbb {R}}^n$ ($n\geq 2$)}, Math.  Scand. \textbf{83} (1998), \mbox{300--308.}\par

\item[{[24]}] N.~Suzuki, \emph{Nonintegrability of harmonic functions in a domain}, Japan. J. Math. \textbf{16} (1990), \mbox{269--278.}\par

\item[{[25]}] 
A.~Torchinsky, \emph{Real-Variable Methods in Harmonic Analysis},  Academic Press, London (1986).\par

\item[{[26]}]
M.~Vuorinen, \emph{On the Harnack constant and the boundary behavior of Harnack functions},  
Ann. Acad. Sci. Fenn., Ser. A I, Math. {\textbf{7}} (1982), 259--277.\par
\end{enumerate}
\end{document}